\documentclass[a4paper,11pt]{amsart}
\usepackage{bm}
\usepackage{amsmath,amssymb}

\usepackage[dvipdfm,bookmarks=true,bookmarksnumbered=true]{hyperref}
\newtheorem{theorem}{Theorem}[section]
\newtheorem{lemma}[theorem]{Lemma}
\newtheorem{proposition}[theorem]{Proposition}
\newtheorem{corollary}[theorem]{Corollary}

\theoremstyle{definition}
\newtheorem{definition}[theorem]{Definition}

\theoremstyle{remark}
\newtheorem{remark}[theorem]{Remark}
\newcommand{\defeq}{\overset{\operatorname{\scriptstyle def.}}{=}}

\newcommand{\beq}{\begin{equation}}
\newcommand{\eeq}{\end{equation}}
\newcommand{\beqa}{\begin{eqnarray}}
\newcommand{\eeqa}{\end{eqnarray}}

\newcommand{\m}{\mu}
\newcommand{\n}{\nu}
\renewcommand{\k}{\kappa}
\newcommand{\lam}{\lambda}
\newcommand{\CP}{{\mathbb P}}
\newcommand{\SU}{\operatorname{\mathrm SU}}

\begin{document}
\pagestyle{headings}
\title{Flop invariance of the topological vertex}
\author{Yukiko Konishi}
\address{Research Institute for Mathematical Sciences, Kyoto University, Kyoto, 606-8502, Japan}
\email{\href{mailto:konishi@kurims.kyoto-u.ac.jp}{konishi@kurims.kyoto-u.ac.jp}}
\author{Satoshi Minabe}
\address{Graduate School of Mathematics, Nagoya University,  Nagoya, 464-8602, Japan}
\email{\href{mailto:minabe@yukawa.kyoto-u.ac.jp}{minabe@yukawa.kyoto-u.ac.jp}}

\subjclass[2000]{Primary 14N35; Secondary 05E05}
\maketitle
\begin{abstract}
We prove transformation formulae for generating functions of Gromov--Witten invariants 
on general toric Calabi--Yau threefolds under flops.  
Our proof is based on a combinatorial identity on the topological vertex 
and analysis of fans of toric Calabi--Yau threefolds.
\end{abstract}
\section{Introduction}
Motivated by a conjecture \cite{Mor, Witten} on quantum cohomology, 
Li and Ruan studied the transformation of Gromov--Witten (GW) invariants 
of projective Calabi--Yau (CY) threefolds under flops 
using symplectic approach \cite{lr}.
The algebro-geometric approach 
was pursued in \cite{ly}. 
The same problem for
Donaldson--Thomas invariants was studied in \cite{HL},
and this may be  related
since
there is a conjecture that 
Donaldson--Thomas invariants and GW invariants 
are related at the level of generating functions \cite{MNOP1, MNOP2}.

In this paper, we study the behavior of GW invariants of 
toric Calabi--Yau (TCY) threefolds (which are noncompact)
under a flop  
based on the method of the topological vertex.
It is a formalism
which expresses the partition functions of 
GW invariants of TCY threefolds 
in terms of symmetric functions \cite{AKMV}.
(In this paper, the partition function 
of GW invariants means the 
exponential of the generating function.) 
Although its original argument was based on the duality to
the Chern-Simons theory,
a mathematical theory including a definition of 
GW invariants for TCY threefolds 
has been developed later in  \cite{LLLZ} 
(see Remark \ref{rem:conj-lllz}). 
We remark that in \cite{IK3}, 
the case of  some special
TCY threefolds was  studied (see remark \ref{rem:ik}).

Let us explain the results of this paper.
Let $X$ be a TCY threefold containing 
a torus invariant rational curve $C$ such that 
its normal bundle is isomorphic to $\mathcal{O}_{\CP^1}(-1) \oplus \mathcal{O}_{\CP^1}(-1)$.
Let $X^+$ be another TCY threefold obtained by flopping $C$. 
Identifying the expansion parameters with respect to second homology classes, 
we can compare the partition function of 
GW invariants of $X$ and that of $X^+$. 
We show that they are equal except for 
factors coming from 
multiples of $[C]$ and from multiples of 
the class $[C^+]$ of the flopped curve $C^+$ (Theorem \ref{prop:flop2}).
Since
the difference between the two appears only 
at the local contributions from neighborhoods of $C$ and $C^+$,
showing the equality of two partition functions
reduces to showing a combinatorial identity 
on skew Schur functions (Theorem \ref{prop:flop}).
Then we obtain the same
result as \cite{lr,ly} on the relation between GW invariants of $X$ and
those of $X^+$.
As an example, we consider the TCY threefold $X$ 
containing two disjoint $\CP^1\times \CP^1$'s
and another related by a flop.
We also show that the partition function of $X$
reproduces Nekrasov's partition function
of $4$-dimensional 
$\SU(2) \times \SU(2)$ gauge theory with a matter
in the bifundamental representation $(\bf{2}, \bar{\bf{2}})$ 
\cite{Nek2} (Proposition \ref{prop:nek}). 
As another application, we consider 
the canonical bundle $K_S$ of a complete smooth toric surface $S$
and the canonical bundle $K_{\hat{S}}$ of 
a blown-up surface $\hat{S}$
and show  that GW invariants of $K_{\hat{S}}$ 
with certain second homology classes 
are equal to those of $K_{S}$ (Proposition \ref{prop:blowup}).

The organization of this paper is as follows. 
In \S \ref{sec:vertex}, we prove a key combinatorial identity.
In \S \ref{sec:tcy}, we give a definition of TCY threefolds used in this paper
and review the method to write down their partition functions.  
In \S \ref{sec:FB}, we study the transformations of partition functions 
under a flop.
In \S \ref{sec:engineering},
we give an example 
and discuss the relationship with 
Nekrasov's partition function. 
In \S \ref{sec:blowup},
we study GW invariants of the canonical bundles of 
smooth toric surfaces related by a blowup.
Combinatorial formulae are collected in Appendix \ref{ap}. 

\subsection*{Acknowledgement}
The authors would like to thank the organizers of the workshop
``Symplectic varieties and related topics'' 
held at Hokkaido University in November 2005. 
This work has grown out from discussions there. 
They are grateful to Prof. H. Kanno 
for useful discussions, comments
and kindly providing them with
his unpublished manuscript \cite{Kanno}.
The result of \S \ref{sec:vertex} is  a generalization
of his result.
They also thank Y. Tachikawa for discussions and comments. 
A part of \S \ref{sec:engineering} is based on discussions with him.

\section{Topological vertex under Flops}\label{sec:vertex}

\subsection{Definitions}
Let $\mathcal{P}$ be the set of partitions. 
For $\mu = (\mu_1 \geq \mu_2 \geq \cdots) \in \mathcal{P}$, 
we define two integers $|\m|$ and $\kappa(\m)$ by
\beq\notag
|\m| = \sum_{j=1}^{l(\m)} \mu_j~, \quad 
\kappa(\m) = |\m| + \sum_{j=1}^{l(\m)} \mu_j (\mu_j -2j)~,
\eeq
where $l(\m)$ is the number of nonzero components in
$\m$.

We use the following definition of the topological vertex 
(\cite{ORV}):
\begin{definition}
\beq
C_{\lam_1, \lam_2, \lam_3}(q) 
\defeq q^{\frac{1}{2}\k(\lam_3)}
s_{\lam_2}(q^\rho) 
\sum_{\m\in\mathcal{P}} s_{\lam_1/\m}(q^{{\lam_2^t} + \rho}) s_{\lam_3^t/ \m} ( q^{{\lam_2} + \rho})~,
\label{topSch}
\eeq
where $s_{\m/\n}(q^{\m + \rho})~\left(\mathrm{resp.}~ s_{\m}(q^{\rho})\right)$ 
is the skew Schur function with the specialization of
variables:
\beq\notag
s_{\m/\n} (x_i = q^{\m_i - i + \frac{1}{2}})~\hskip4mm 
(\mathrm{resp.} ~s_{\m} (x_i = q^{- i + \frac{1}{2}}))~.
\eeq
\end{definition}

Take four partitions $\lambda_1, \lambda_2, \lambda_3, \lambda_4$.
These will be fixed throughout the rest of \S \ref{sec:vertex}.
We define
\begin{eqnarray}
Z_{0}(q,Q_0)&\defeq&\sum_{\mu\in\mathcal{P}}(-Q_{0})^{|\mu|}
C_{\lambda_1, \lambda_2, \mu^{t}}(q)
C_{\lambda_{3}, \lambda_{4}, \mu}(q)~, \\
Z_{0}^{+}(q,Q_0^+)&\defeq&\sum_{\mu\in\mathcal{P}}(-Q_{0}^{+})^{|\mu|}
C_{\lambda_1,  \mu^t, \lambda_4}(q)
C_{\lambda_{3}, \mu, \lambda_2}(q)~.
\end{eqnarray}
We also set
\begin{equation}\label{eq:conifold-partition-fcn}
Z_{(-1,-1)}(q,Q)
  =\prod_{k=1}^{\infty}(1-Q q^k)^k,
\end{equation}
and
\begin{equation*}
Z_{0}^{\prime}(q,Q_0) \defeq 
\frac{Z_0(q,Q_0)}{Z_{(-1,-1)}(q,Q_0)}~, ~\hskip4mm 
Z_{0}^{+\prime}(q,Q_0^+) \defeq 
\frac{Z_0^+(q,Q_0^+)}{Z_{(-1,-1)}(q,Q_0^+)}~.
\end{equation*}
The goal of this section is to show  an identity 
relating  $Z_{0}'(q,Q_0)$ and $Z_{0}^{+\prime}(q,Q_0^+)$ 
under the identification $Q_0^+=Q_0^{-1}$ (Theorem \ref{prop:flop}).  
Formulae necessary for  proofs can be found in Appendix \ref{ap}.

\begin{remark}
Let us mention the geometrical meaning of the above formal power series.
$Z_{(-1,-1)}(q,Q_0)$ is the partition function 
of the TCY threefold 
$\mathcal{O}_{\CP^1}(-1)\oplus \mathcal{O}_{\CP^1}(-1)$
(cf. \S \ref{sec:partition-fcn},
see also \cite[(C.18)]{EK1} and \cite[Theorem 3]{FP}).
$Z_{0}(q,Q_0)$ 
and $Z_0^+(q,Q_0^+)$ appear as local contributions in
the partition functions of TCY threefolds
related by a flop 
such that both a flopping curve and a 
flopped curve
have  normal bundles 
isomorphic to
$\mathcal{O}_{\CP^1}(-1)\oplus \mathcal{O}_{\CP^1}(-1)$
(see  Figure \ref{fig:flop-fan}).
\end{remark}
\subsection{Individual calculations}
First, we compute $Z_{0}'(q, Q_0)$ and $Z_{0}^{+\prime}(q, Q_0^+)$ 
respectively.  

Let us introduce the following functions: 
\beqa\notag
&&f_{\m}(q)={q \over q-1}\sum_{i\geq1}(q^{\mu_i-i}-q^{-i})~,\\
\notag
&&f_{\m, \n}(q)=(q-2+q^{-1})f_{\m}(q)f_{\n}(q)+f_{\m}(q)+f_{\n}(q)~,
\eeqa
and let $C_k(\m,\n)$ be  the expansion coefficients in
the Laurent polynomial ${f}_{\m,\n}(q)$:
\beq\notag
f_{\m, \n}(q)
=\sum_{k\in \mathbb{Z}} C_k(\m, \n) q^k~.
\eeq \label{prod}

\begin{proposition}
We have 
\begin{equation}
\begin{split}\label{eq:Z-before-flop}
Z_{0}'(q,Q_0)
&=
q^{\frac{1}{2}\kappa(\lambda_2)+\frac{1}{2}\kappa(\lambda_4)}
s_{\lambda_1}(q^{\rho})s_{\lambda_3}(q^{\rho}) 
\prod_{k\in\mathbb{Z}}(1-Q_{0} q^k)^{C_k(\lambda_1^t,\lambda_3^t)}\\
&\times 
\sum_{\tau}(-Q_0)^{|\tau|}
s_{\lambda_2^t/\tau^t}(q^{\lambda_1+\rho}, Q_0 q^{-\lambda_3-\rho})
s_{\lambda_4^t/\tau}(q^{\lambda_3+\rho}, Q_0 q^{-\lambda_1-\rho})~,
\end{split}
\end{equation}
\begin{equation}
\begin{split}
Z_{0}^{+\prime}(q,Q_0^+)
&=
s_{\lambda_1}(q^{\rho}) s_{\lambda_3}(q^{\rho}) 
\prod_{k\in\mathbb{Z}}(1-Q_{0}^{+} q^k)^{C_k(\lambda_1, \lambda_3)}\\
&\times
\sum_{\tau}(-Q_{0}^{+})^{|\tau|}
s_{\lambda_2/\tau}(q^{\lambda_3^t+\rho}, Q_{0}^{+} q^{-\lambda_1^t-\rho})
s_{\lambda_4/\tau^t}(q^{\lambda_1^t+\rho}, Q_{0}^{+} q^{-\lambda_3^t-\rho})~.
\end{split}
\label{eq:Zflop}
\end{equation} 
\label{prop:z}
\end{proposition}

\begin{proof}
By definition (\ref{topSch}) of the topological vertex, we have
\begin{eqnarray*}
\begin{split}
Z_{0}(q,Q_0)&=\sum_{\mu}(-Q_0)^{|\mu|}
q^{\frac{1}{2}\kappa(\lambda_2)}s_{\lambda_1}(q^{\rho})
\sum_{T}s_{\mu^t/T}(q^{\lambda_1^t+\rho})s_{\lambda_2^t/T}(q^{\lambda_1+\rho})\\
&\qquad \qquad  q^{\frac{1}{2}\kappa(\lambda_4)}s_{\lambda_3}(q^{\rho})
\sum_{T'}s_{\mu/T'}(q^{\lambda_3^t+\rho})s_{\lambda_4^t/T'}(q^{\lambda_3+\rho})\\
&=q^{\frac{1}{2}\kappa(\lambda_2)+\frac{1}{2}\kappa(\lambda_4)}
s_{\lambda_1}(q^{\rho})s_{\lambda_3}(q^{\rho})
\sum_{T,T'}(-Q_0)^{|T|}
s_{\lambda_2^t/T}(q^{\lambda_1+\rho})s_{\lambda_4^t/T'}(q^{\lambda_3+\rho})\\
&\qquad \qquad \sum_{\mu}s_{\mu^t/T}(-Q_0 q^{\lambda_1^t+\rho})s_{\mu/T'}(q^{\lambda_3^t+\rho})~.
\end{split}
\end{eqnarray*}
We perform the sum with respect to $\m$ by using  (\ref{Schur2}):
\begin{eqnarray*}
\begin{split}
Z_{0}(q,Q_0)
&=q^{\frac{1}{2}\kappa(\lambda_2)+\frac{1}{2}\kappa(\lambda_4)}
       s_{\lambda_1}(q^{\rho})s_{\lambda_3}(q^{\rho})\\
&\qquad \prod_{i,j \geq 1}(1-Q_{0} q^{h_{\lambda_1^t, \lambda_3^t}(i,j)})
\sum_{\tau}s_{T^t/\tau}(q^{\lambda_3^t+\rho})s_{(T')^t/\tau^t}(-Q_0 q^{\lambda_1^t+\rho})\\
&\qquad \sum_{T,T'}(-Q_0)^{|T|}
s_{\lambda_2^t/T}(q^{\lambda_1+\rho})s_{\lambda_4^t/T'}(q^{\lambda_3+\rho})\\
&=q^{\frac{1}{2}\kappa(\lambda_2)+\frac{1}{2}\kappa(\lambda_4)}
   s_{\lambda_1}(q^{\rho})s_{\lambda_3}(q^{\rho})
   \prod_{i,j \geq 1}(1-Q_{0} q^{h_{\lambda_1^t, \lambda_3^t}(i,j)})\\
&\qquad \sum_{\tau} (-Q_0)^{|\tau|}
\sum_{T}s_{\lambda_2^t/T}(q^{\lambda_1+\rho})s_{T/\tau^t}(Q_0 q^{-\lambda_3-\rho})
\sum_{T'}s_{\lambda_4^t/T'}(q^{\lambda_3+\rho})s_{T'/\tau}(Q_0 q^{-\lambda_1-\rho})~.
\end{split}
\end{eqnarray*}
Here for $\mu,\nu\in\mathcal{P}$,
\beq\notag
h_{\m, \n}(i,j) \defeq \m_i -i  +\n_j - j +1~.
\label{relative hook}\eeq
In passing to the second line, we have used (\ref{symmetry}). 
By using (\ref{Schur3}), we have 
\begin{eqnarray*}
\begin{split}
Z_{0}(q,Q_0)=q^{\frac{1}{2}\kappa(\lambda_2)+\frac{1}{2}\kappa(\lambda_4)}
&s_{\lambda_1}(q^{\rho})s_{\lambda_3}(q^{\rho}) 
\prod_{i,j \geq 1}(1-Q_{0} q^{h_{\lambda_1^t, \lambda_3^t}(i,j)})\\
&\sum_{\tau}(-Q_0)^{|\tau|}
s_{\lambda_2^t/\tau^t}(q^{\lambda_1+\rho}, Q_0 q^{-\lambda_3-\rho})
s_{\lambda_4^t/\tau}(q^{\lambda_3+\rho}, Q_0 q^{-\lambda_1-\rho})~.
\end{split}
\end{eqnarray*}
Applying Lemma \ref{lem:h-C}, we obtain (\ref{eq:Z-before-flop}). 
One can also compute $Z_{0}^+(q,Q_0^+)$ in a similar way.
\end{proof}

The next corollary is a consequence of
Proposition \ref{prop:z}.
\begin{corollary}
$Z_{0}^{+\prime}(q,Q_0^+)$ is a polynomial in $Q_0^+$
 of degree  at most
$|\lambda_1|+|\lambda_2|+|\lambda_3|+|\lambda_4|$. 
Moreover,
if $\lambda_3=\lambda_4=\emptyset$,
$Z_{0}^{+\prime}(q,Q_0^+)$ is a polynomial in $Q_0^+$
of degree $|\lambda_1|+|\lambda_2|$. 
\label{lem:polynomial}
\end{corollary}
Similar statement also holds  for $Z_0'(q,Q_0)$.

\begin{proof}
The first statement follows if we apply (\ref{eq:sumc}) and (\ref{homogeneity})
to the expression
(\ref{eq:Zflop}).
To prove the second statement, we show that
the top term does not vanish. 
By (\ref{eq:sumc}), we have
\begin{equation*}
\prod(1-Q_0^+q^k)^{C_k(\lambda_1^t, \emptyset)} 
= (-1)^{|\lambda|} q^{-\frac{1}{2}\kappa_{\lambda_1}}(Q_0^+)^{|\lambda_1|} 
+ (\text{terms of  lower degree in $Q_0^+$}).
\end{equation*} 
Substituting this into  (\ref{eq:Zflop}) 
with $\lambda_3,\lambda_4$ set to $\emptyset$,
and using (\ref{homogeneity}),
we obtain the claim.
\end{proof}

\subsection{Comparison}
Next, we compare $Z_{0}'(q,Q_0)$ with $Z_{0}^{+\prime}(q,Q_0^+)$ 
under the identification  $Q_{0}^{+}=Q_0^{-1}$. 
First we have the following

\begin{lemma}
Under the identification $Q_{0}^{+}=Q_0^{-1}$, we have
\begin{multline}
\sum_{\tau}(-Q_{0}^{+})^{|\tau|}s_{\lambda_2/\tau}(q^{\lambda_3^t+\rho}, Q_{0}^{+} q^{-\lambda_1^t-\rho})
s_{\lambda_4/\tau^t}(q^{\lambda_1^t+\rho}, Q_{0}^{+} q^{-\lambda_3^t-\rho})\\
=(-Q_0)^{-|\lambda_2|-|\lambda_4|} \sum_{\tau}(-Q_0)^{|\tau|}
s_{\lambda_2^t/\tau^t}(q^{\lambda_1+\rho}, Q_0 q^{-\lambda_3-\rho})
s_{\lambda_4^t/\tau}(q^{\lambda_3+\rho}, Q_0 q^{-\lambda_1-\rho})~.
\label{eq:sumpart}
\end{multline}
\label{lem:sum}
\end{lemma}

\begin{proof}
Under  $Q_{0}^{+}=Q_0^{-1}$, we have 
\begin{eqnarray*}
(\mathrm{LHS})
&=&\sum_{\tau}(-Q_0^{-1})^{|\tau|}
s_{\lambda_2/\tau}(q^{\lambda_3^t+\rho}, Q_0^{-1} q^{-\lambda_1^t-\rho})
s_{\lambda_4/\tau^t}(q^{\lambda_1^t+\rho}, Q_0^{-1} q^{-\lambda_3^t-\rho})\\
&=&\sum_{\tau}(-Q_0^{-1})^{|\tau|}(-1)^{|\lambda_2|+|\lambda_4|}
s_{\lambda_2^t/\tau^t}(q^{-\lambda_3-\rho}, Q_0^{-1} q^{\lambda_1+\rho})
s_{\lambda_4^t/\tau}(q^{-\lambda_1-\rho}, Q_0^{-1} q^{\lambda_3+\rho})~, \\
&=&\sum_{\tau}(-Q_0)^{-|\lambda_2|+|\tau|-|\lambda_4|+|\tau|} (-Q_0^{-1})^{|\tau|}
s_{\lambda_2^t/\tau^t}(q^{\lambda_1+\rho}, Q_0 q^{-\lambda_3^t-\rho})
s_{\lambda_4^t/\tau}(q^{\lambda_3+\rho}, Q_0 q^{-\lambda_1^t-\rho})\\
&=& (\mathrm{RHS})~.
\end{eqnarray*}
Note that we have used the property (\ref{symmetry}) in the second line and
the homogeneity (\ref{homogeneity}) of skew Schur functions in the third line.
\end{proof}

The next lemma was proven in \cite[eq.(45)]{IK3}.
\begin{lemma}
The following identity holds:
\begin{eqnarray*}
\prod_{k}(1-Q_0^{-1}q^k)^{C_k(\lambda_1, \lambda_3)}
=(-Q_0)^{-|\lambda_1|-|\lambda_3|} q^{\frac{1}{2}\kappa(\lambda_1)+\frac{1}{2}\kappa(\lambda_3)}
\prod_{k}(1-Q_0q^k)^{C_k(\lambda_1^t, \lambda_3^t)}.
\end{eqnarray*}
\label{lem:product}
\end{lemma}

\begin{proof}
By (\ref{eq:c}), we have 
\begin{eqnarray*}
\prod_{k}(1-Q_0^{-1}q^k)^{C_k(\lambda_1, \lambda_3)}
&=&\prod_{k}(1-Q_0^{-1}q^{-k})^{C_k(\lambda_1^t, \lambda_3^t)}\\
&=& Q_0^{-\frac{1}{2}(|\lambda_1|+|\lambda_3| )}
 q^{-\frac{1}{4}(\kappa(\lambda_1^t)+\kappa(\lambda_3^t))}
\prod_{k}\left( Q_0^{\frac{1}{2}}q^{\frac{k}{2}} - Q_0^{-\frac{1}{2}}q^{-\frac{k}{2}}\right)
^{C_k(\lambda_1^t, \lambda_3^t)}\, .
\end{eqnarray*}
On the other hand, we have 
\begin{eqnarray*}
\prod_{k}(1-Q_0q^{k})^{C_k(\lambda_1^t, \lambda_3^t)}
&=&Q_0^{\frac{1}{2}(|\lambda_1|+|\lambda_3| )}
 q^{\frac{1}{4}(\kappa(\lambda_1^t)+\kappa(\lambda_3^t))}
\prod_{k}\left( Q_0^{-\frac{1}{2}}q^{-\frac{k}{2}} - Q_0^{\frac{1}{2}}q^{\frac{k}{2}}\right)
^{C_k(\lambda_1^t, \lambda_3^t)}\\
&=&(-1)^{|\lambda_1|+|\lambda_3|} Q_0^{\frac{1}{2}(|\lambda_1|+|\lambda_3| )}
 q^{\frac{1}{4}(\kappa(\lambda_1^t)+\kappa(\lambda_3^t))}
\prod_{k}\left( Q_0^{\frac{1}{2}}q^{\frac{k}{2}} - Q_0^{-\frac{1}{2}}q^{-\frac{k}{2}}\right)
^{C_k(\lambda_1^t, \lambda_3^t)}\, .
\end{eqnarray*}
By comparing the above two equations and 
by using a symmetry (\ref{kappa}) of a $\kappa$-factor, we get the claim. 
\end{proof}

The following is the main result in this section.  
(The case of $\lambda_1=\lambda_4=\emptyset$ was proved in 
\cite{Kanno}.) 
\begin{theorem}
Under the identification $Q_{0}^{+} =Q_0^{-1}$, we have
\begin{eqnarray*}
Z_{0}^{+\prime}(q,Q_0^+) = 
(-Q_0)^{-(|\lambda_1|+|\lambda_2|+|\lambda_3|+|\lambda_4|)}
q^{\frac{1}{2}\left(\kappa(\lambda_1)-\kappa(\lambda_2)+\kappa(\lambda_3)-\kappa(\lambda_4)\right)}
Z_{0}^{\prime}(q,Q_0)~ .
\end{eqnarray*}
\label{prop:flop}
\end{theorem}
\begin{proof}
This follows from Proposition \ref{prop:z} and Lemmas \ref{lem:sum} and \ref{lem:product}.
\end{proof}

\section{Toric Calabi--Yau threefolds and Partition Functions}
\label{sec:tcy}
In this section, we give definitions of toric Calabi--Yau threefolds and 
their partition functions. 
Our reference is  \cite{ko2}. 
\subsection{Toric Calabi--Yau threefolds} 
\begin{definition}\label{def:TCY}
A toric Calabi--Yau (TCY) threefold is 
a three-dimensional smooth toric variety  $X$
over $\mathbb{C}$
associated with  a fan $\Sigma$ satisfying following conditions:
\begin{enumerate}
\item[(i)]
the primitive generator $\vec{\omega}$ of every 1-cone
satisfies
$\vec{\omega}\cdot \vec{u}=1$ where $\vec{u}=(0,0,1)$;
\item[(ii)]
all maximal cones are three dimensional;
\item[(iii)]
$|\Sigma|\cap \{z=1\}$ 
is simply connected where 
$|\Sigma|=\displaystyle{\bigcup_{\sigma\in \Sigma}}
\sigma \subset \mathbb{R}^3$ is 
the support of $\Sigma$ and 
$z$ is the third coordinate of $\mathbb{R}^3$.
\end{enumerate}
\end{definition}
The condition (i) is equivalent to the condition that
$\wedge^3 T^*X$ is trivial (Calabi--Yau condition)
and the condition (ii) implies that $\pi_1(X)=0$.
The condition (iii) is imposed for simplicity of arguments.

We briefly describe necessary facts on (co)homology
of TCY threefolds.
Recall that 
the subset $\Sigma_n\subset\Sigma$ 
of $n$-cones is in one-to-one correspondence
with the set of $(3-n)$-dimensional torus invariant subvarieties in $X$.
Let $\Sigma_1=\{\rho_1,\ldots,\rho_r\}$
be the set of 1-cones.
Denote by $\vec{\omega}_i$ $(1\leq i\leq r)$
the primitive lattice vector generating $\rho_i$
and by $D_{\rho_i}\subset X$ ($1\leq i\leq r$)
the torus invariant Weil divisor 
corresponding to $\rho_i$.
The group $A_2(X)$ of all Weil divisors modulo 
rational equivalence
is generated by $D_{\rho_1},\ldots,D_{\rho_r}$ with rational equivalence 
given by
$\displaystyle{
\sum_{j=1}^r A_{ij}D_{\rho_j}=0
}$
$(i=1,2,3)$
(\cite{fulton}, the first proposition in \S 3.4)
where $A=(A_{ij})$ is the $3\times r$ matrix
\begin{equation}\notag
A=\big(
\vec{\omega}_1,\ldots,\vec{\omega}_r\big).
\end{equation}
Let $\Sigma_2'$ be the set of 
2-cones which lie in the interior of $|\Sigma|$:
\begin{equation}\notag
\Sigma_2'=\{\tau\in \Sigma_2|\, 
\tau\subset |\Sigma|\setminus\partial |\Sigma|\}.
\end{equation}
It is in one-to-one correspondence with the set of 
torus invariant (hence rational) curves in $X$.
Let us write $\Sigma_2'=\{\tau_1,\ldots,\tau_p\}$
and let $C_{\tau_i}\subset X$ denote 
the rational curve corresponding to
$\tau_i$.
We define $N_1^T(X)$ to
be the set of 2-cycles generated by $C_{\tau_1},\ldots,C_{\tau_p}$
modulo numerical equivalence.
Note that by the intersection pairing
$A_2(X)\times N_1^T(X) \to \mathbb{Z}
$,
$A_2(X)\otimes \mathbb{R}$ and
$N_1^T(X)\otimes \mathbb{R}$ become  dual to each other.

Now let us explain the calculation of the intersection numbers
and numerical equivalence.
If $\rho_j$ and  $\tau_i$ spans a 3-cone,
$D_{\rho_j}.C_{\tau_i}=1$  and
if
$\rho_j$ and  $\tau_i$ do not span a cone in the fan,
$D_{\rho_j}.C_{\tau_i}=0$
(\cite{fulton},\S 5,1 p.98).
If two 1-cones,  say $\rho_1,\rho_2$, are contained in $\tau_i$,
then
$D_{\rho_1}.C_{\tau_i}$ and $D_{\rho_2}.C_{\tau_i}$ are  obtained
via rational equivalence relations of $D_{\rho_j}$'s.
For convenience, we  introduce
the following injective map
\begin{equation}
l_X:N_1^T(X)\to 
    \{l\in \mathbb{Z}^{r}\,| \, A.l=\vec{0}\}=L_A
,\quad
Z\mapsto(D_{\rho_1}.Z,\ldots,D_{\rho_r}.Z).
\end{equation}
Then
$D_{\rho_1}.C_{\tau_i}$ and $D_{\rho_2}.C_{\tau_i}$ are  obtained 
by solving the equation $A.l_X([C_{\tau_i}])=\vec{0}$.
(Hence they satisfy the relation  
$D_{\rho_1}.C_{\tau_i}+D_{\rho_2}.C_{\tau_i}=-2$.)
The numerical equivalence
can be read  from  linear relations between the vectors 
$l_X([C_{\tau_1}]),\dots,l_X([C_{\tau_p}])$.

By the analysis of  the gluing of local coordinate systems around 
$C_{\tau_i}$, we see that
its normal bundle is isomorphic to
$\mathcal{O}_{\mathbb{P}^1}
(D_{\rho_1}.C_{\tau_i})\oplus
\mathcal{O}_{\mathbb{P}^1}(D_{\rho_2}.C_{\tau_i})$.
We will use a term
a $(-1,-1)$-curve
for a torus invariant curve with the 
normal bundle isomorphic to
$\mathcal{O}_{\mathbb{P}^1}
(-1)\oplus
\mathcal{O}_{\mathbb{P}^1}(-1)$.

\subsection{Partition functions}\label{sec:partition-fcn}

Let  $X$ be a TCY threefold and $\Sigma$ be its fan.
We briefly review how to write down the 
partition function of $X$.

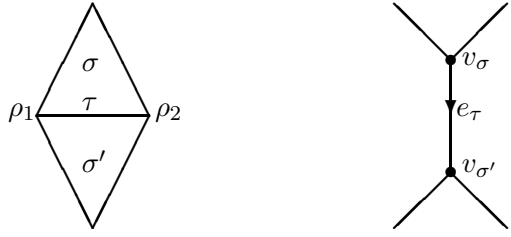
\begin{figure}[t]
\begin{center}
\unitlength 0.15cm
\begin{picture}(10,20)(0,-10)
\thicklines
\put(0,0){\line(1,0){10}}
\put(0,0){\line(1,2){5}}
\put(10,0){\line(-1,2){5}}
\put(0,0){\line(1,-2){5}}
\put(10,0){\line(-1,-2){5}}
\put(-2.5,0){$\rho_1$}
\put(10.5,0){$\rho_2$}
\put(4,.5){$\tau$}
\put(4,4){$\sigma$}
\put(4,-5){$\sigma'$}
\end{picture}
\hspace*{3cm}
\begin{picture}(10,20)(0,-10)
\thicklines
\put(5,-5){\line(0,1){5}}
\put(5,-5){\line(-1,-1){5}}
\put(5,-5){\line(1,-1){5}}
\put(5,5){\vector(0,-1){5}}
\put(5,5){\line(-1,1){5}}
\put(5,5){\line(1,1){5}}
\put(5,-5){\circle*{1}}
\put(5,5){\circle*{1}}
\put(5.5,0){$e_{\tau}$}
\put(6,4.5){$v_{\sigma}$}
\put(6,-5){$v_{\sigma'}$}
\end{picture}
\end{center}
\caption{Fan (section at $z=1$) and toric graph}
\label{fig:fan-graph}
\end{figure}

First,
consider the following directed graph $\Gamma_X$ (called  a toric graph)
with labels on edges of a certain type.
The vertex set is
\begin{equation}\notag
V(\Gamma_X)=V_3(\Gamma_X)\cup V_1(\Gamma_X),
\,
V_3(\Gamma_X)=\{v_{\sigma}|\sigma\in \Sigma_3(X)\},
\,
V_1(\Gamma_X)=\{v_{\tau}|\tau\in \Sigma_2(X)\setminus\Sigma_2'(X)\}.
\end{equation}
The edge set is
\begin{equation}\notag
E(\Gamma_X)=E_3(\Gamma_X)\cup E_1(\Gamma_X),
\,
E_3(\Gamma_X)=\{e_{\tau}|\tau\in \Sigma_2'(X)\},
\,
E_1(\Gamma_X)=
\{e_{\tau}|\tau\in \Sigma_2(X)\setminus \Sigma_2'(X)\}.
\end{equation}
An edge $e_{\tau}\in E_3(\Gamma_X)$ joins
$v_{\sigma},v_{\sigma'}\in V_3(\Gamma)$ iff
$\tau=\sigma\cap\sigma'$ (see Figure \ref{fig:fan-graph})
and an edge
$e_{\tau}\in E_1(\Gamma)$ joins
$v_{\sigma}\in V_3(\Gamma_X)$ and $v_{\tau}\in V_1(\Gamma_X)$
iff $\sigma$ is a unique 3-cone
such that $\tau$ is a face of $\sigma$.
(Note that a vertex in $V_3(\Gamma_X)$ is trivalent
and a vertex in $V_1(\Gamma_X)$ is univalent.)
The direction of edges can be taken arbitrarily.
The label $n:E_3(\Gamma)\to \mathbb{Z}$,
called the {\em framing}, is given as follows:
\begin{equation}\notag
n(e_{\tau})=\frac{D_{\rho_1}.C_{\tau}-D_{\rho_2}.C_{\tau}}{2}~,
\end{equation}
where
$\tau\in\Sigma_2'$ and
$\rho_1,\rho_2\in \Sigma_1$ are as shown in Figure \ref{fig:fan-graph}.
Note that $\Gamma_X$ is connected by the condition (iii) in
Definition \ref{def:TCY}

\newcommand{\threein}{  
\unitlength .15cm
\begin{picture}(5,12)(-5,-5)
\thicklines
\put(0,0){\line(1,1){5}}
\put(0,0){\line(-1,1){5}}
\put(0,0){\line(0,-1){5}}
\put(0,0){\vector(1,1){4}}
\put(0,0){\vector(-1,1){4}}
\put(0,0){\vector(0,-1){4}}
\put(0,0){\circle*{1}}
\put(4,2){$e'$}
\put(-5.5,1.5){$e^{\prime\prime}$}
\put(1,-4){$e$}
\put(-2,-1){$v$}
\end{picture}
}
\newcommand{\twoin}{  
\unitlength .15cm
\begin{picture}(5,5)(-5,-5)
\thicklines
\put(0,0){\line(1,1){5}}
\put(0,0){\line(-1,1){5}}
\put(0,0){\line(0,-1){5}}
\put(0,0){\vector(1,1){4}}
\put(0,0){\vector(-1,1){4}}
\put(0,-5){\vector(0,1){3}}
\put(4,2){$e'$}
\put(-5.5,1.5){$e^{\prime\prime}$}
\put(1,-4){$e$}
\put(0,0){\circle*{1}}
\put(-2,-1){$v$}
\end{picture}
}
\newcommand{\onein}{  
\unitlength .15cm
\begin{picture}(5,5)(-5,-5)
\thicklines
\put(0,0){\line(1,1){5}}
\put(0,0){\line(-1,1){5}}
\put(0,0){\line(0,-1){5}}
\put(0,0){\vector(-1,1){4}}
\put(5,5){\vector(-1,-1){3}}
\put(0,-5){\vector(0,1){3}}
\put(0,0){\circle*{1}}
\put(4,2){$e'$}
\put(-5.5,1.5){$e^{\prime\prime}$}
\put(1,-4){$e$}
\put(-2,-1){$v$}
\end{picture}
}
\newcommand{\threeout}{  
\unitlength .15cm
\begin{picture}(5,5)(-5,-5)
\thicklines
\put(0,0){\line(1,1){5}}
\put(0,0){\line(-1,1){5}}
\put(0,0){\line(0,-1){5}}
\put(5,5){\vector(-1,-1){3}}
\put(-5,5){\vector(1,-1){3}}
\put(0,-5){\vector(0,1){3}}
\put(0,0){\circle*{1}}
\put(4,2){$e'$}
\put(-5.5,1.5){$e^{\prime\prime}$}
\put(1,-4){$e$}
\put(-2,-1){$v$}
\end{picture}
}
Secondly, we write down the partition function
from $\Gamma_X$.
Let 
\begin{equation}\notag
\mathcal{P}(\Gamma_X)=\{\vec{\lambda}: E_3(\Gamma)\to \mathcal{P}\}.
\end{equation}
Take the set of formal variables
$\vec{Q}=(Q_e)_{e\in E_3(\Gamma_X)}$ associated to
$E_3(\Gamma_X)$.
Then the partition function of $X$
is a formal power series in $\vec{Q}$ given by
\begin{equation}\label{eq:partition-fct}
Z_X(q,\vec{Q})
=\sum_{\vec{\lambda}\in \mathcal{P}(\Gamma)}
\prod_{e\in E_3(\Gamma)}(-1)^{|\vec{\lambda}(e)|(n_e+1)}
q^{\frac{\kappa(\vec{\lambda}(e))}{2}n(e)}
Q_e^{|\vec{\lambda}(e)|}
\prod_{v\in V_3(\Gamma)}C_{\vec{\lambda}_v}(q)~.
\end{equation}
Here
$C_{\vec{\lambda}_v}(q)$ is the topological vertex 
defined in (\ref{topSch})
and
$\Vec{\lambda}_v$ ($v\in V_3(\Gamma)$,
$\Vec{\lambda}\in \mathcal{P}(\Gamma)$)
is as in Figure \ref{fig:lambdav}
(for $e\in E(\Gamma_X)\setminus E_3(\Gamma_X)$, set $\vec{\lambda}(e)$
to $\emptyset$).
We remark that the partition function does not depend on
the directions of edges
since  the framing changes the sign
if one gives the opposite direction to an edge $e\in E_3(\Gamma_X)$
and it is compensated by (\ref{kappa}) and the summation.

\begin{figure}[t]
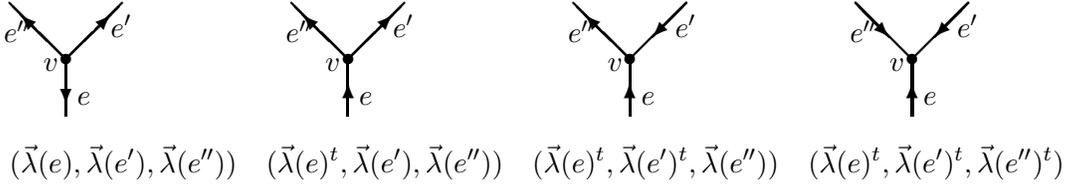

\begin{equation}\notag
\begin{split}
&
\threein
\hspace*{3cm}
\twoin
\hspace*{3cm}
\onein
\hspace*{3cm}
\threeout
\\
&(\vec{\lambda}(e),\Vec{\lambda}(e'),\vec{\lambda}(e^{\prime\prime}))\quad
(\vec{\lambda}(e)^t,\Vec{\lambda}(e'),\vec{\lambda}(e^{\prime\prime}))\quad
(\vec{\lambda}(e)^t,\Vec{\lambda}(e')^t,\vec{\lambda}(e^{\prime\prime}))\quad
(\vec{\lambda}(e)^t,\Vec{\lambda}(e')^t,\vec{\lambda}(e^{\prime\prime})^t)
\end{split}
\end{equation}
\caption{$\vec{\lambda}_v$}\label{fig:lambdav}
\end{figure}

\begin{remark}\label{rem:conj-lllz}
Precisely speaking, the partition function obtained in 
\cite{LLLZ} has the expression 
almost same as (\ref{eq:partition-fct}) 
except that  $C_{\vec{\lambda}_v}(q)$ is replaced by
$\tilde{\mathcal{W}}_{\vec{\lambda}_v}(q)$. 
Here $\tilde{\mathcal{W}}_{\lambda_1, \lambda_2,\lambda_3}(q)$
is a rational function in $q^{\frac{1}{2}}$
similar to $C_{\lambda_1, \lambda_2,\lambda_3}(q)$
but has a slightly different expression. 
It is conjectured that 
$\tilde{\mathcal{W}}_{\lambda_1, \lambda_2,\lambda_3}(q) =C_{\lambda_1, \lambda_2,\lambda_3}(q)$
\cite[Conjecture 8.3]{LLLZ}. Here we use $C_{\lambda_1, \lambda_2,\lambda_3}(q)$
assuming that the conjecture is true. 
\end{remark}

The Gromov--Witten invariant $N_{g,\beta}(X)$ of $X$ with the genus $g$
and the second homology class $\beta\in H_2^{cpt}(X,\mathbb{Z})$ 
(see \cite{LLLZ} for a definition)
is obtained as follows:
\begin{equation}\label{eq:GW}
\sum_{g\geq 0} N_{g,\beta}(X) g_s^{2g-2}=
\sum_{
\begin{subarray}{c}\vec{d}=(d_e)_{e\in E_3(\Gamma_X)},\\ 
                   \vec{d} [\vec{C}]=[\beta]
\end{subarray}}
{F}_{\vec{d}}(e^{\sqrt{-1}g_s}),
\end{equation}
where $[\vec{C}]=([C_e])_{e\in E_3(\Gamma_X)}$
and $C_e\subset X$ is the rational curve  corresponding to $e$.
${F}_{\vec{d}}(q)$
is the coefficient of 
$\vec{Q}^{\vec{d}}=\displaystyle{\prod_{e\in E_3(\Gamma_X)}}
Q_e^{d_e}$
in $\log Z_X(q,\vec{Q})$.

\section{Transformations of partition functions under flop} 
\label{sec:FB}
In this section, we study the transformation of the partition function
of  TCY threefolds under a flop. 

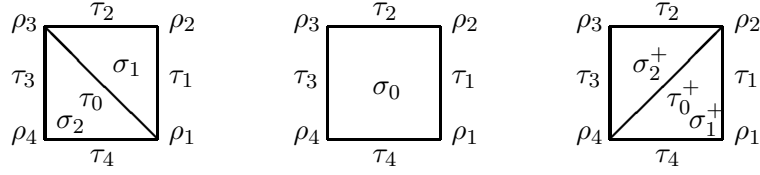
\begin{figure}

\vspace*{5mm}
\begin{center}
\unitlength .15cm

\begin{picture}(10,10)
\thicklines
\put(0,0){\line(1,0){10}}
\put(0,0){\line(0,1){10}}
\put(10,10){\line(0,-1){10}}
\put(10,10){\line(-1,0){10}}
\put(10,0){\line(-1,1){10}}
\put(3,3){$\tau_0$}
\put(11,5){$\tau_1$}
\put(4,11){$\tau_2$}
\put(-3,5){$\tau_3$}
\put(4,-2){$\tau_4$}
\put(11,0){$\rho_1$}
\put(11,10){$\rho_2$}
\put(-3,10){$\rho_3$}
\put(-3,0){$\rho_4$}
\put(1,1){$\sigma_2$}
\put(6,6){$\sigma_1$}
\end{picture}
\hspace*{2cm}
\begin{picture}(10,10)
\thicklines
\put(0,0){\line(1,0){10}}
\put(0,0){\line(0,1){10}}
\put(10,10){\line(0,-1){10}}
\put(10,10){\line(-1,0){10}}
\put(11,5){$\tau_1$}
\put(4,11){$\tau_2$}
\put(-3,5){$\tau_3$}
\put(4,-2){$\tau_4$}
\put(11,0){$\rho_1$}
\put(11,10){$\rho_2$}
\put(-3,10){$\rho_3$}
\put(-3,0){$\rho_4$}
\put(4,4){$\sigma_0$}
\end{picture}
\hspace*{2cm}
\begin{picture}(10,10)
\thicklines
\put(0,0){\line(1,0){10}}
\put(0,0){\line(0,1){10}}
\put(10,10){\line(0,-1){10}}
\put(10,10){\line(-1,0){10}}
\put(0,0){\line(1,1){10}}
\put(5,3){$\tau_0^+$}
\put(11,5){$\tau_1$}
\put(4,11){$\tau_2$}
\put(-3,5){$\tau_3$}
\put(4,-2){$\tau_4$}
\put(11,0){$\rho_1$}
\put(11,10){$\rho_2$}
\put(-3,10){$\rho_3$}
\put(-3,0){$\rho_4$}
\put(2,6){$\sigma_2^+$}
\put(7,1){$\sigma_1^+$}
\end{picture}

\end{center}
\vspace*{5mm}

\caption{Fans (sections at $z=1$): $\Sigma$ (left),
$\bar{\Sigma}$ (middle) and $\Sigma^+$ (right).
The generators $\vec{\omega}_1,\ldots,\vec{\omega}_4$
of $\rho_1,\dots,\rho_4$ satisfy the relation 
$\vec{\omega}_1+\vec{\omega}_3=\vec{\omega}_2+\vec{\omega}_4$.}
\label{fig:flop-fan}
\end{figure}

Let $X$ be a TCY threefold and let $\Sigma$ be its fan.
Assume that $X$ contains at least one $(-1,-1)$-curve $C_0$.
Denote the corresponding 2-cone by $\tau_0$.
Near $\tau_0$, the fan looks like 
the left diagram in Figure \ref{fig:flop-fan}.
We set
\begin{equation}\notag
\bar{\Sigma}=(\Sigma\setminus\{\tau_0,\sigma_1,\sigma_2\})
\cup\{\sigma_0\},
\quad
\Sigma^+=
(\Sigma\setminus\{\tau_0,\sigma_1,\sigma_2\})
\cup \{\tau^+_0,\sigma_1^+,\sigma_2^+\}
\end{equation}
where $\tau_0,\sigma_1,\sigma_2,\sigma_0,
\tau_0^+\sigma_1^+,\sigma_2^+$ are cones shown in
Figure \ref{fig:flop-fan}.
Let $Y$ be
the singular toric variety associated with the fan $\bar{\Sigma}$
and $X^+$ be the TCY threefold associated with the fan $\Sigma^+$.
Then 
associated to the evident maps
$\Sigma\to \bar{\Sigma}$ and
$\Sigma^+ \to \bar{\Sigma}^+$,
there are the following birational maps:
\begin{equation}\notag
\begin{matrix}
X&         &\stackrel{\phi}{\dashrightarrow}&             &\!\!\!\!X^+\\
 &\!\!\!\!\!f\searrow &&\swarrow f^+\!\!\!\!\!&   \\
 &       & Y &       &
\end{matrix}
\end{equation}
The map $f$ is a  small contraction with the exceptional set $C_0$
and $\phi$ is a flop of $f$.
\begin{remark}\label{rmk:phi*}
Since $\Sigma$ and $\Sigma^+$ have the same set of 1-cones,
there is a canonical isomorphism $A_2(X)\cong A_2(X^+)$
induced by $\phi$.
In turn, this induces an isomorphism 
$\phi_*: N_1^T(X)\otimes \mathbb{R}\to
N_1^T(X^+)\otimes \mathbb{R}$ via the duality between
$A_2(\cdot)\otimes \mathbb{R}$ and 
$N_1^T(\cdot)\otimes \mathbb{R}$ where $\cdot=X,X^+$.
\end{remark}

From here on,
we proceeds assuming that $\tau_1,\ldots,\tau_4\in \Sigma_2'$.
Other cases can be recovered by setting to zero the formal variables
associated to any of $\tau_1,\dots,\tau_4$ which are not in $\Sigma_2'$.
We use the notations shown in Table \ref{tab:notation-flop}.

\begin{table}
\begin{equation}\notag
\begin{array}{|c|cc|cc|}\hline
& X && X^+&\\\hline
\text{2-cone}
   & \tau_0,\tau_1,\dots,\tau_4  & \tau 
   & \tau_0^+,\tau_1,\dots,\tau_4& \tau  
\\
\text{curve}
   & C_0,C_1,\dots,C_4  & C_{\tau}
   & C_0^+,C_1^+,\dots,C_4^+ & C_{\tau}^+             
\\
\text{edge} 
      & e_0,e_1,\dots,e_4 & e_{\tau}\text{ or just }e 
      & e_0^+,e_1^+,\dots,e_4^+ &e_{\tau}\text{ or just }e 
\\
\text{variable}
      & Q_0,Q_1,\dots,Q_4 & Q_e
      & Q_0^+,Q_1^+,\dots,Q_4^+ & Q_e
\\\hline
\end{array} 
\end{equation}

\caption{}
\label{tab:notation-flop}
\end{table}

\begin{lemma}\label{lem:curves-flop}
Under the flop
$\phi:X \dashrightarrow  X^+$,
the curve classes transform as follows.
\begin{equation}\notag
\phi_*[C_0]=-[C^+_0],\qquad
\phi_*[C_i]=[C_i^+]+[C^+_0],
\qquad
\phi_*[C_{\tau}]=[C_{\tau}^+] \quad
(\tau\in \Sigma_2'(X)\setminus\{\tau_0,\dots,\tau_4\}).
\end{equation}
\end{lemma}
\begin{proof}
The first statement follows from
$l_X([C_0])=-l_{X^+}([C_0^+])$
by  remark \ref{rmk:phi*}.
The proof of the other two is similar.
\end{proof}

\begin{figure}[t]
\begin{center}
\unitlength .15cm
\begin{picture}(30,30)(-10,-10)
\thicklines

\put(10,10){\line(1,0){10}}
\put(10,10){\vector(1,0){6}}
\put(14,8){$e_{1}$}

\put(0,0){\line(0,-1){10}}
\put(0,0){\vector(0,-1){6}}
\put(1,-6){$e_4$}

\put(0,0){\line(-1,0){10}}
\put(0,0){\vector(-1,0){6}}
\put(-7,1){$e_3$}

\put(10,10){\line(0,1){10}}
\put(10,10){\vector(0,1){6}}
\put(11,15){$e_2$}

\put(0,0){\line(1,1){10}}
\put(0,0){\vector(1,1){6}}
\put(7,4){$n(e_0)=0$}
\put(2,6){$e_{0}$}

\put(0,0){\circle*{1}}
\put(10,10){\circle*{1}}

\put(-2,1){$v_2$}
\put(7,10){$v_1$}
\end{picture}
\hspace*{2cm}
%
\begin{picture}(40,30)(-10,-10)
\thicklines

\put(10,0){\line(1,0){10}}
\put(10,0){\vector(1,0){6}}
\put(15,1){$e_{1}^+$}
\put(12,-3){$n(e_{1}^+)=n(e_{1})-1$}

\put(10,0){\line(0,-1){10}}
\put(10,0){\vector(0,-1){6}}
\put(6,-8){$e_{4}^+$}
\put(12,-8){$n(e_{4}^+)=n(e_{4})+1$}

\put(0,10){\line(-1,0){10}}
\put(0,10){\vector(-1,0){6}}
\put(-9,11.5){$e_{3}^+$}
\put(1,12){$n(e_{3}^+)=n(e_{3})-1$}

\put(0,10){\line(0,1){10}}
\put(0,10){\vector(0,1){6}}
\put(-4,17){$e_{2}^+$}
\put(1,17){$n(e_{2}^+)=n(e_{2})+1$}

\put(10,0){\line(-1,1){10}}
\put(0,10){\vector(1,-1){6}}
\put(6,6){$n(e_{0}^+)=0$}
\put(2,2.5){$e_{0}^+$}

\put(10,0){\circle*{1}}
\put(0,10){\circle*{1}}

\put(6.5,-2){$v_1^+$}
\put(-2.5,7.5){$v_2^+$}

\end{picture}
\end{center}

\caption{Toric graphs $\Gamma_X$ (left)
and $\Gamma_{X^+}$ (right).}\label{TF1}
\end{figure}
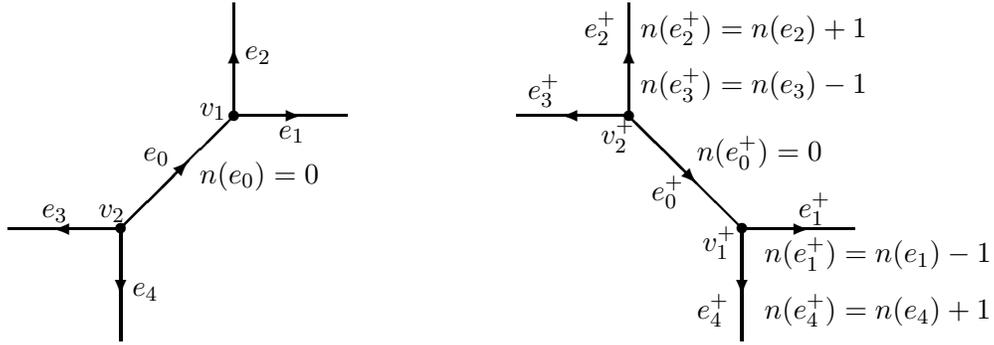

Let $\Gamma_X$ be a toric graph of $X$. 
Near the edge $e_{0}$, the graph looks like 
the left diagram in Figure \ref{TF1}.
Under the flop $\phi$, the toric diagram (and
the framings) changes as follows.
\begin{lemma}
A graph obtained from $\Gamma_X$ by replacing
the left diagram in Figure \ref{TF1}
with the right is a toric graph of $X^+$.
\end{lemma}

We associate the same formal variables 
$\vec{Q}=(Q_e)$
to edges in
$E_3(\Gamma_X)\setminus \{e_0,\dots,e_4\}$
and
those in $E_3(\Gamma_{X^+})\setminus \{e_0^+,\dots,e_4^+\}$
and write
the partition functions of $X$ and $X^+$
as 
$Z_X(q,\vec{Q},Q_{0},Q_1,Q_2,Q_3,Q_4)$
and
$Z_{X^+}(q,\vec{Q},Q_0^+,Q_1^+,Q_2^+,Q_3^+,Q_4^+)$ respectively.
It is immediate to check that
\begin{equation}
\begin{split}\label{floppingcurve}
Z_X(q,\vec{0},Q_0,0,0,0,0)&=
Z_{(-1,-1)}(q,Q_0),\\
Z_{X^+}(q,\vec{0},Q_0^+,0,0,0,0)
&=
Z_{(-1,-1)}(q,Q_0^+).
\end{split}
\end{equation}
We set 
\begin{eqnarray*}
Z_X^{\prime}(q,\vec{Q},Q_{0},Q_1,Q_2,Q_3,Q_4)
&\defeq&
\frac{Z_X(q,\vec{Q},Q_{0},Q_1,Q_2,Q_3,Q_4)}{Z_X(q,\vec{0},Q_0,0,0,0,0)}~,\\
Z_{X^+}^{\prime}(q,\vec{Q},Q_0^+,Q_1^+,Q_2^+,Q_3^+,Q_4^+)
&\defeq&
\frac{Z_{X^+}(q,\vec{Q},Q_0^+,Q_1^+,Q_2^+,Q_3^+,Q_4^+)}
{Z_{X^+}(q,\vec{0},Q_0^+,0,0,0,0)}~.
\end{eqnarray*}

Now we will compare these.
To do so,
we should identify the formal variables
so that the identification is compatible with Lemma \ref{lem:curves-flop}:
\begin{equation}\notag
Q_{0}=(Q_{0}^+)^{-1},\qquad Q_{i}=Q_{0}^+ Q_{i}^+.
\end{equation}

\begin{theorem}

\begin{enumerate}
\item[(i)]
The coefficients of $\vec{Q}^{\vec{d}}Q_0^{d_0}Q_1^{d_1}\dots Q_4^{d_4}$
in $Z_X^{\prime}(q,\vec{Q},Q_{0},Q_1,Q_2,Q_3,Q_4)$
is zero if $d_0>d_1+d_2+d_3+d_4$.
A similar result holds for $X^+$.

\item[(ii)]
Under the identification
$Q_{0}=(Q_{0}^+)^{-1}$,
$Q_{i}=Q_{0}^+ Q_{i}^+$, we have
\begin{equation}\notag
Z_X^{\prime}(q,\vec{Q},Q_{0},Q_1,Q_2,Q_3,Q_4)
=
Z_{X^+}^{\prime}(q,\vec{Q},Q_0^+,Q_1^+,Q_2^+,Q_3^+,Q_4^+)~.
\end{equation}
(This is an equality between
two formal power series in $\vec{Q},Q_0^+,\dots,Q_4^+$.)
\label{prop:flop2}
\end{enumerate}

\end{theorem}

\begin{remark}\label{rem:ik}
In \cite[\S 4.1]{IK3}, 
Iqbal and Kashani--Poor studied the 
special case such that 
the 2-cones $\tau_2,\tau_4\notin\Sigma_2'$ and 
the curves $C_1,C_3$ have normal bundles
$\mathcal{O}_{\CP^1}(-1)\oplus \mathcal{O}_{\CP}(-1)$
or $\mathcal{O}_{\CP^1}(-2)\oplus \mathcal{O}_{\CP^1}(0)$.
They obtained the result of Lemma \ref{lem:curves-flop} 
and proved the second statement of Theorem \ref{prop:flop2}
in that case.
\end{remark}

\begin{proof}

(i) follows from the first statement of Corollary \ref{lem:polynomial}.

(ii)
Let
\begin{equation}\notag
\mathcal{P}'(\Gamma_X)=
 \{\vec{\nu}:  E_3(\Gamma_X)\setminus\{e_0\}\to \mathcal{P}\}
\end{equation}
and define
$\vec{\nu}_v\in \mathcal{P}^3$ 
for $\vec{\nu}\in \mathcal{P}'(\Gamma_X)$ and 
$v\in V_3(\Gamma_X)\setminus\{v_1,v_2\} $
in the same way as 
$\vec{\lambda}_v$ (Figure \ref{fig:lambdav}).
After (\ref{eq:partition-fct}),
$Z_X(q,\vec{Q},Q_{0},Q_1,Q_2,Q_3,Q_4)$ 
is written as follows:
\begin{equation}\notag
\begin{split}
&Z_X(q,\vec{Q},Q_{0},Q_1,Q_2,Q_3,Q_4)\\
&=\sum_{\vec{\nu}\in \mathcal{P}'(\Gamma_X)}
\prod_{e\in E_3(\Gamma_X)\setminus\{e_0,\dots,e_4\}}
(-1)^{(n(e)+1)|\vec{\nu}(e)|}Q_e^{|\vec{\nu}(e)|}
\prod_{v\in V_3(\Gamma_X)\setminus\{v_1,v_2\}}C_{\vec{\nu}_v}(q)
\\
&\times
\underbrace{ 
\prod_{i=1}^4 (-1)^{(n(e_i)+1)|\vec{\nu}(e_i)|}Q_i^{|\vec{\nu}(e_i)|}
\sum_{\mu\in \mathcal{P}}
          C_{\vec{\nu}(e_1),\vec{\nu}(e_2),\mu^t}(q)
          C_{\vec{\nu}(e_3),\vec{\nu}(e_4),\mu}(q) (-Q_0)^{|\mu|}
}_{(a)}.
\end{split}
\end{equation}

Similarly,
$Z_{X^+}(q,\vec{Q},Q_{0}^+,Q_1^+,Q_2^+,Q_3^+,Q_4^+)$ 
is written as follows:
\begin{equation}\notag
\begin{split}
&Z_{X^+}(q,\vec{Q},Q_{0}^+,Q_1^+,Q_2^+,Q_3^+,Q_4^+)\\
&=\sum_{\vec{\nu}\in \mathcal{P}'(\Gamma_{X^+})}
\prod_{e\in E_3(\Gamma_{X^+})\setminus\{e_0^+,\dots,e_4^+\}}
(-1)^{(n(e)+1)|\vec{\nu}(e)|}Q_e^{|\vec{\nu}(e)|}
\prod_{v\in V_3(\Gamma_{X^+})\setminus\{v_1^+,v_2^+\}}C_{\vec{\nu}_v}(q)
\\
&\times
\underbrace{ 
\prod_{i=1}^4 (-1)^{(n(e_i^+)+1)|\vec{\nu}(e_i^+)|}
(Q_i^+)^{|\vec{\nu}(e_i^+)|}
\sum_{\mu\in \mathcal{P}}
          C_{\vec{\nu}(e_1^+),\mu^t,\vec{\nu}(e_4^+)}(q)
          C_{\vec{\nu}(e_3^+),\mu,\vec{\nu}(e_2^+)}(q) (-Q_0^+)^{|\mu|}
}_{(b)}.
\end{split}
\end{equation}
Here 
\begin{equation}\notag
\mathcal{P}'(\Gamma_{X^+})=
 \{\vec{\nu}:  E_3(\Gamma_{X^+})\setminus\{e_0^+\}\to \mathcal{P}\}
\end{equation}
and for $\vec{\nu}\in \mathcal{P}'(\Gamma_{X^+})$ and 
$v\in V_3(\Gamma_{X^+})\setminus\{v_1^+,v_2^+\} $,
$\vec{\nu}_v\in \mathcal{P}^3$ is defined in the same way.

Since $\Gamma_X$ and $\Gamma_{X^+}$ are identical outside the 
diagrams described in Figure \ref{TF1},
$E_3(\Gamma_{X})\setminus\{e_0,\dots,e_4\}=
 E_3(\Gamma_{X^+})\setminus\{e_0^+,\dots,e_4^+\}
$, 
$ V_3(\Gamma_{X})\setminus\{v_1,v_2\}=
 V_3(\Gamma_{X^+})\setminus\{v_1^+,v_2^+\}
$
and we have a natural bijection
$p:\mathcal{P}'(\Gamma_{X})\to\mathcal{P}'(\Gamma_{X^+})$
such that $p(\vec{\nu})=\vec{\nu}^+$ iff $\vec{\nu}(e)=\vec{\nu}^+(e)$
for all $e\in E_3(\Gamma_{X})\setminus\{e_0,\dots,e_4\}$
and $\vec{\nu}(e_i)=\vec{\nu}^+(e_i^+)$ for $1\leq i\leq 4$.
Under this identification,
we could see that the two partition functions
have the same expressions except for the factors $(a)$ and $(b)$.
Taking into account the change in framings,
we have
\begin{equation}\notag
\frac{(a)}{Z_{(-1,-1)}(q,Q_0)}\,\Bigg|_{Q_0=(Q_0^+)^{-1},Q_i=Q_0^+Q_i^+}
=
\frac{(b)}{Z_{(-1,-1)}(q,Q_0^+)}
\end{equation}
by Theorem \ref{prop:flop}.
\end{proof}

We finish this subsection by
restating Theorem \ref{prop:flop2} in terms of GW invariants.
(Compare with \cite[Corollary A.1]{lr} and \cite[Theorem 3.1.1]{ly}.)
\begin{corollary}\label{cor:GW-flop}
For $\beta\in H_2^{cpt}(X,\mathbb{Z})$
such that $\beta$ is not a multiple of $[C_0]$,
\begin{equation}\notag
N_{g,\phi_*(\beta)}(X^+)=N_{g,\beta}(X).
\end{equation}
Moreover,
\begin{equation}\notag
N_{g,d[C_0]}(X)=N_{g,d[C_0^+]}(X^+)
=N_{g,d[\CP^1]}(\mathcal{O}_{\CP^1}(-1)\oplus\mathcal{O}_{\CP^1}(-1)).
\end{equation}
\end{corollary}

\begin{proof}
Theorem \ref{prop:flop2} implies that
$\log Z_X(q,\vec{Q},Q_0,\dots,Q_4)$ and
$\log Z_{X^+}(q,\vec{Q},Q_0^+,\dots,Q_4^+)$
are  written in the following form:
\begin{equation}\notag
\begin{split}
\log Z_X'(q,\vec{Q},Q_0,\dots,Q_4)
&=\sum_{\vec{d}}
\sum_{\begin{subarray}{c}d_0,\ldots,d_4\geq 0,\\
d_1+\dots+d_4\geq d_0\end{subarray}} F_{\vec{d},d_0,d_1,d_2,d_3,d_4}(q)
\vec{Q}^{\vec{d}}
Q_0^{d_0}\dots Q_4^{d_4},
\\
\log Z_{X^+}'(q,\vec{Q},Q_0^+,\dots,Q_4^+)
&=\sum_{\vec{d}}
\sum_{\begin{subarray}{c}d_0,\ldots,d_4\geq 0,\\
d_1+\dots+d_4\geq d_0\end{subarray}} F_{\vec{d},d_0,d_1,d_2,d_3,d_4}(q)
\vec{Q}^{\vec{d}}
(Q_0^+)^{d_1+\dots+d_4-d_0}
(Q_1^+)^{d_1}\dots (Q_4^+)^{d_4}.
\end{split}
\end{equation}
Comparing with (\ref{eq:GW}), we obtain the first statement.
The second statement follows from (\ref{floppingcurve}).
\end{proof}

\section{Example and geometric engineering} \label{sec:engineering}

In this section, we first give an example of
\S \ref{sec:FB}.
Then  we will discuss its relation with
Nekrasov's partition function \cite{Nek2}
along the same lines with \cite{IK1, IK2, EK1, EK2, Zhou}.

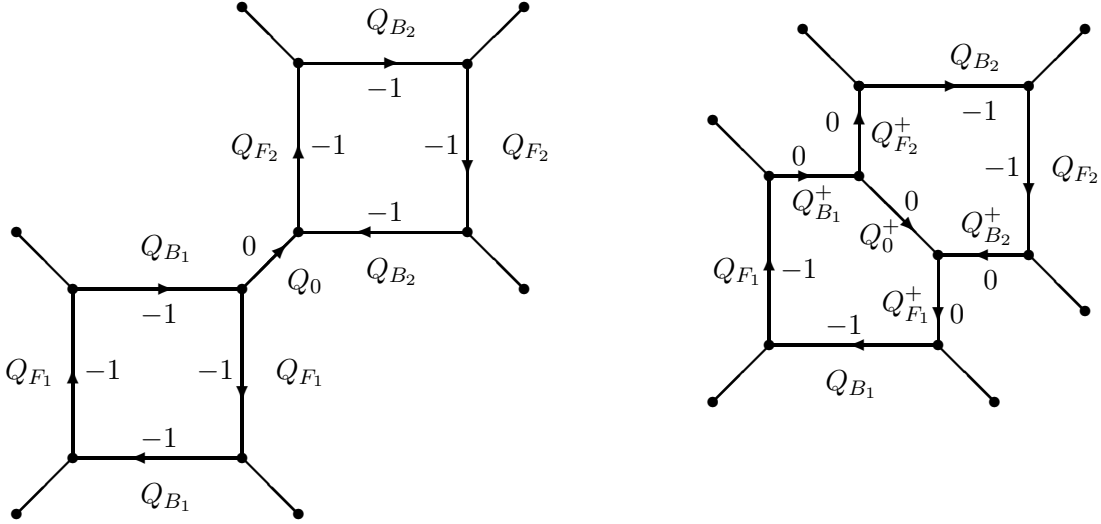
\begin{figure}[t]
\begin{center}
\unitlength .15cm
\begin{picture}(30,30)(-10,-10)
\thicklines

\put(-5,-5){\line(-1,1){5}}
\put(-5,-20){\line(1,0){15}}
\put(10,-20){\line(0,1){15}}
\put(5, -20){\vector(-1,0){5}}

\put(-5,-5){\circle*{1}}
\put(-10,0){\circle*{1}}
\put(-5,-5){\line(1,0){15}}
\put(-2,-5){\vector(1,0){6}}

\put(-5,-20){\circle*{1}}
\put(-10,-25){\circle*{1}}
\put(-5,-20){\line(0,1){15}}
\put(-5,-20){\line(-1,-1){5}}
\put(-5,-15){\vector(0,1){3}}

\put(10,-20){\circle*{1}}
\put(15,-25){\circle*{1}}
\put(10,-20){\line(1,-1){5}}
\put(10, -5){\vector(0,-1){10}}

\put(10,-5){\circle*{1}}
\put(10,-5){\line(1,1){5}}
\put(11,-4){\vector(1,1){3}}
\put(15,15){\line(-1,1){5}}
\put(15,0){\line(1,0){15}}
\put(30,0){\line(0,1){15}}
\put(25, 0){\vector(-1,0){5}}

\put(30,15){\circle*{1}}
\put(35,20){\circle*{1}}
\put(15,15){\line(1,0){15}}
\put(18,15){\vector(1,0){6}}

\put(15,0){\circle*{1}}
\put(10,20){\circle*{1}}
\put(15,0){\line(0,1){15}}
\put(15,0){\line(-1,-1){5}}
\put(15,5){\vector(0,1){3}}

\put(30,0){\circle*{1}}
\put(35,-5){\circle*{1}}
\put(30,0){\line(1,-1){5}}
\put(30, 15){\vector(0,-1){10}}

\put(15,15){\circle*{1}}
\put(30,15){\line(1,1){5}}

\put(1,-24){$Q_{B_1}$}
\put(1,-19){$-1$}
\put(-11,-13){$Q_{F_1}$}
\put(-4,-13){$-1$}
\put(1, -2){$Q_{B_1}$}
\put(1,-8){$-1$}
\put(13,-13){$Q_{F_1}$}
\put(6,-13){$-1$}
\put(14, -5){$Q_0$}
\put(10,-2){$0$}
\put(21,-4){$Q_{B_2}$}
\put(21,1){$-1$}
\put(9,7){$Q_{F_2}$}
\put(16,7){$-1$}
\put(21, 18){$Q_{B_2}$}
\put(21,12){$-1$}
\put(33,7){$Q_{F_2}$}
\put(26,7){$-1$}

\end{picture}
\hspace{3cm}
%
\begin{picture}(40,0)(-10,-10)
\thicklines
\put(5,5){\line(-1,1){5}}
\put(5,-10){\line(1,0){15}}

\put(5,5){\circle*{1}}
\put(0,10){\circle*{1}}
\put(6,5){\vector(1,0){3}}

\put(5,5){\line(1,0){8}}
\put(13,5){\circle*{1}}
\put(0,-15){\circle*{1}}
\put(25,-15){\circle*{1}}
\put(33,-7){\circle*{1}}
\put(33,18){\circle*{1}}
\put(8,18){\circle*{1}}
\put(13,7){\vector(0,1){4}}
\put(14,4){\vector(1,-1){4}}

\put(5,-10){\circle*{1}}
\put(5,-10){\line(0,1){15}}
\put(5,-10){\line(-1,-1){5}}
\put(5,-5){\vector(0,1){3}}

\put(20,-10){\circle*{1}}
\put(20,-10){\line(1,-1){5}}
\put(17,-10){\vector(-1,0){5}}
\put(20,-5){\vector(0,-1){3}}
\put(13,13){\line(-1,1){5}}

\put(28,-2){\line(0,1){15}}

\put(13,13){\circle*{1}}
\put(13,13){\circle*{1}}
\put(13,13){\line(1,0){15}}
\put(16,13){\vector(1,0){6}}

\put(13,5){\line(0,1){8}}

\put(28,-2){\circle*{1}}
\put(28,-2){\line(1,-1){5}}
\put(28, 13){\vector(0,-1){10}}
\put(28,-2){\line(-1, 0){8}}
\put(28,13){\circle*{1}}
\put(28,13){\line(1,1){5}}
\put(27, -2){\vector(-1,0){4}}
\put(20,-2){\circle*{1}}
\put(20,-2){\line(0,-1){8}}
\put(20,-2){\line(-1,1){7}}

\put(10,-14){$Q_{B_1}$}
\put(10,-9){$-1$}
\put(7,2){$Q_{B_1}^+$}
\put(7,6){$0$}
\put(0, -4){$Q_{F_1}$}
\put(6,-4){$-1$}
\put(15,-7){$Q_{F_1}^+$}
\put(21,-8){$0$}

\put(13, -1){$Q_0^+$}
\put(17,2){$0$}

\put(22,0){$Q_{B_2}^+$}
\put(24,-5){$0$}
\put(14, 8){$Q_{F_2}^+$}
\put(10,9){$0$}
\put(21, 15){$Q_{B_2}$}
\put(22,10){$-1$}
\put(30,5){$Q_{F_2}$}
\put(24, 5){$-1$}
\end{picture}

\end{center}
\vspace{3cm}
\caption{TCY threefold which contains two disjoint $\CP^1 \times \CP^1$ connected by a $(-1,-1)$-curve (left) and its flop (right).}
\label{fig:quiver}
\end{figure}

Let $X$ and $X^+$ be 
the TCY threefolds associated with 
the left and right toric graphs in Figure \ref{fig:quiver}, respectively. 
$X$ 
contains two copies of $\CP^1\times \CP^1$ disjoint to each other
and $X^+$ is  obtained  by a flop of a unique $(-1,-1)$-curve in $X$.
In this example,
formal variables should be assigned as in Figure \ref{fig:quiver}:
the five variables for $X$ are independent
and the nine variables for $X^+$ have the four relations
$Q_{F_i}=Q_{F_i}^+Q_0^+$ ($i=1,2$)
and $Q_{B_i}=Q_{B_i}^+Q_0^+$ ($i=1,2$).
The variables of $X$ and $X^+$
should be identified by
$Q_0^+=Q_0^{-1}$ and $Q_{F_i},Q_{B_i}$ of $X^+$
$=Q_{F_i},Q_{B_i}$ of $X$.

Let us compute the partition function 
$Z_X$ of $X$ (we omit the variables).
By Proposition \ref{prop:z},  we have  
\begin{eqnarray*}
\begin{split}
&Z_X  
=\sum_{\mu^1_1, \mu^1_2,\mu^2_1,\mu^2_2}~
\prod_{k=1}^2~
 (Q_{B_k})^{|\mu^k_1|+|\mu^k_2|}
s_{\mu^k_1}^2(q^{\rho}) s_{\mu^k_2}^2(q^{\rho}) 
\prod_{i,j \geq1}\left(1-Q_{F_k} q^{h_{\mu^k_1, (\mu^k_2)^t}(i,j)}\right)^{-2}\\
&\qquad \qquad
\sum_{\lambda} ~(-Q_0)^{|\lambda|} ~
s_{\lambda^t}(q^{\mu^1_2+\rho},  Q_{F_1}q^{\mu^1_1+\rho})~
s_{\lambda}(q^{(\mu^2_2)^t +\rho},  Q_{F_2}q^{(\mu^2_1)^t +\rho})~.
\end{split}
\end{eqnarray*}
We can perform the sum 
in the last factor by (\ref{Schur2}):  
\begin{eqnarray*}
\sum_{\lambda} &(-Q_0)^{|\lambda|} 
s_{\lambda^t}(q^{\mu^1_2+\rho},  Q_{F_1}q^{\mu^1_1+\rho})
s_{\lambda}(q^{(\mu^2_2)^t +\rho},  Q_{F_2}q^{(\mu^2_1)^t +\rho})&\\
= ~\prod_{i,j\geq1}
&\left( 1- Q_0 q^{h_{\mu^1_2, (\mu^2_2)^t}(i,j)}\right)
\left( 1- Q_0 Q_{F_1} q^{h_{\mu^1_1, (\mu^2_2)^t}(i,j)}\right)&\\
&\left( 1- Q_0 Q_{F_2} q^{h_{\mu^1_2, (\mu^2_1)^t}(i,j)}\right)
\left(1- Q_0 Q_{F_1} Q_{F_2} q^{h_{\mu^1_1, (\mu^2_1)^t}(i,j)}\right)&~.
\end{eqnarray*} 
Therefore
Theorem \ref{prop:flop2} implies that
the partition function $Z_{X^+}$ is obtained from $Z_X$ by
replacing
\begin{equation}\notag
\prod_{i,j\geq1} \left( 1- Q_0 q^{h_{\mu^1_2, (\mu^2_2)^t}(i,j)}\right)
\, \to \, 
\prod_k \left( 1- (Q_0^+)^{-1} q^k \right)^{C_k(\mu^1_2, (\mu^2_2)^t)}
\prod_{k\geq 1}\left( 1- Q_0^+ q^{k}\right)^k,
\end{equation}
and replacing $Q_0$ in other factors by $(Q_0^+)^{-1}$.

From the discussions in \cite[\S 2.1]{KMV}, 
it seems natural to expect that  
the partition function of $X$ reproduces 
Nekrasov's partition function for a gauge theory with a product gauge group
and with a matter.
We want to clarify this statement.  
Let us set 
\begin{equation*}
Z_{X}^{\mathrm{inst}}=\frac{Z_X}{Z_X|_{Q_{B_1} = Q_{B_2} = 0}~.}
\end{equation*}
Then, by the same method with \cite{EK1, Zhou}, we can show the following
\begin{proposition}
Let 
\begin{equation*}\label{eq:sub}
q=e^{-2R\hbar}~,\quad
Q_{F_1}=e^{-4R a_{1}}~,\quad
Q_{F_2}=e^{-4R a_{2}}~,\quad
Q_0 = e^{2R(a_{1}+a_{2}-m)}~.
\end{equation*}
Then we have
\begin{eqnarray*}
\begin{split}
Z_X^{\mathrm{inst}} 
&= \quad
\sum_{\mu^1_1, \mu^1_2, \mu^2_1, \mu^2_2}~ 
\prod_{k=1}^2 \left( \frac{Q_{B_k}}{2^4 Q_{F_k}} \right)^{|\mu^k_1| + |\mu^k_2|}
\prod_{l,n =1}^{2} \prod_{i,j\geq 1}~
\frac{\sinh R\left( a^{(k)}_{ln} + \hbar \left( \mu^k_{1,i} - \mu^k_{2,j} +j-i \right)  \right)}
{\sinh R\left( a^{(k)}_{ln}  + \hbar \left( j-i \right)  \right)} \\
&
\times q^{\frac{1}{2} 
\left( \kappa(\mu^1_1) +\kappa(\mu^1_2) - \kappa(\mu^2_1) -\kappa(\mu^2_2) \right)}~
(2^2 Q_0)^{|\mu^1_1| + |\mu^1_2| + |\mu^2_1|+|\mu^2_2|}~
(Q_{F_1}^{\frac{1}{2}})^{ 2|\mu^1_1|+|\mu^2_1|+|\mu^2_2| }~
(Q_{F_2}^{\frac{1}{2}})^{|\mu^1_1|+|\mu^1_2|+2|\mu^2_1|}\\
& \qquad \qquad \qquad \qquad \qquad \qquad \qquad
\prod_{l,n=1}^2~ \prod_{i,j\geq1}~
\frac{\sinh R\left( a^{(1,2)}_{ln} + m + \hbar \left( j-i \right) \right)}
{\sinh R\left ( a^{(1,2)}_{ln} +m + \hbar \left( \mu^1_{l,i} - \mu^2_{n, j}+j-i \right) \right)}~,
\end{split}
\end{eqnarray*}
where 
\begin{eqnarray*}
a^{(k)}_{11} = a^{(k)}_{22} = 0~, \hskip5mm a^{(k)}_{12} = - a^{(k)}_{21} = 2a_k~,
\end{eqnarray*}
and
\begin{eqnarray*}
a^{(1,2)}_{11}=a_1 +a_2~, \hskip5mm
a^{(1,2)}_{21}=-a_1+a_2~,\hskip5mm
a^{(1,2)}_{12}=a_1-a_2~,\hskip5mm
a^{(1,2)}_{22}=-a_1-a_2~. 
\end{eqnarray*}
\label{prop:nek}
\end{proposition}

By Proposition \ref{prop:nek}, 
it is easy to see that the $R\to0$ limit of
\begin{equation*} 
Z_X^{\mathrm{inst}} 
|_{q=e^{-2R\hbar},~Q_{B_k} =2^2 \Lambda_k,~Q_{F_k}=e^{-4R a_{k}},~Q_0 = e^{2R(a_{1}+a_{2}-m)}}
\end{equation*}
is equal to the instanton part of 
Nekrasov's partition function of $4$-dimensional 
$\SU(2) \times \SU(2)$ gauge theory with a matter
in the bifundamental representation $(\bf{2}, \bar{\bf{2}})$ 
\cite[(66)]{Nek2}. 
(See also \cite{FMP, HIV, Matsuura-Ohta, S} for related works.)

\begin{remark}
It is immediate to see that $Z_{X^+}^{\mathrm{inst}}=
Z_{X^+}/(Z_{X^+}|_{Q_{B_1}=Q_{B_2}=0})$
also coincides with the same Nekrasov's partition function
with a similar variable identification in the limit $R\to 0$.
More generally,
Theorem \ref{prop:flop2} may imply that
if TCY threefolds $X$ and $X^+$ are 
related by flops with respect to $(-1,-1)$-curves
and if the partition function of $X$ reproduces 
Nekrasov's partition function for a gauge theory,
then the partition function of $X^+$ also reproduces it.
(This statement itself seems to be well-known to specialists.)
\end{remark}

\section{Application to toric surface and its blowup}\label{sec:blowup}

\begin{figure}

\vspace*{5mm}
\begin{center}
\unitlength .15cm

\begin{picture}(10,10)
\thicklines
\put(10,10){\line(0,-1){10}}
\put(10,10){\line(-1,0){10}}
\thinlines
\put(10,0){\line(-1,1){10}}
\thicklines
\put(3,3){$\tau_0$}
\put(11,5){$\tau_1$}
\put(4,11){$\tau_2$}
\put(11,0){$\rho_1$}
\put(11,10){$\rho_2$}
\put(-3,10){$\rho_3$}
\put(6,6){$\sigma$}
\end{picture}
\hspace*{2cm}
\begin{picture}(10,10)
\thinlines
\put(0,0){\line(1,0){10}}
\put(0,0){\line(0,1){10}}
\thicklines
\put(10,10){\line(0,-1){10}}
\put(10,10){\line(-1,0){10}}
\put(10,0){\line(-1,1){10}}
\put(3,3){$\tau_0$}
\put(11,5){$\tau_1$}
\put(4,11){$\tau_2$}
\put(-3,5){$\tau_3$}
\put(4,-2){$\tau_4$}
\put(11,0){$\rho_1$}
\put(11,10){$\rho_2$}
\put(-3,10){$\rho_3$}
\put(-3,0){$\rho_4$}
\put(1,1){$\sigma_2$}
\put(6,6){$\sigma_1$}
\end{picture}
\hspace*{2cm}
\begin{picture}(10,10)
\thicklines
\thinlines
\put(0,0){\line(1,0){10}}
\put(0,0){\line(0,1){10}}
\thicklines
\put(10,10){\line(0,-1){10}}
\put(10,10){\line(-1,0){10}}
\put(0,0){\line(1,1){10}}
\put(5,3){$\hat{\tau}_0$}
\put(11,5){$\tau_1$}
\put(4,11){$\tau_2$}
\put(-3,5){$\tau_3$}
\put(4,-2){$\tau_4$}
\put(11,0){$\rho_1$}
\put(11,10){$\rho_2$}
\put(-3,10){$\rho_3$}
\put(-3,0){$\rho_4$}
\put(2,6){$\hat{\sigma}_2$}
\put(7,1){$\hat{\sigma}_1$}
\end{picture}

\end{center}

\vspace*{5mm}

\caption{Fans (sections at $z=1$): $\Sigma$ (left),
$\bar{\Sigma}$ (middle)
and $\hat{\Sigma}$ (right).
The generators $\vec{\omega}_1,\ldots,\vec{\omega}_4$
of $\rho_1,\dots,\rho_4$ satisfy the relation 
$\vec{\omega}_1+\vec{\omega}_3=\vec{\omega}_2+\vec{\omega}_4$.
}
\label{fig:blowup-fan}

\end{figure}
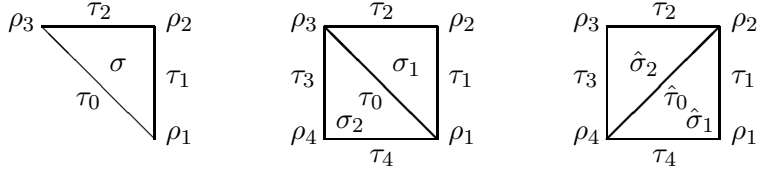
As an application, we compare
GW invariants of the canonical bundle
of a complete smooth toric surface
and
those of the canonical bundle of a blown-up surface.
Some  relevant numerical data 
can be found in \cite{CKYZ}. 

Let $S$ be a complete smooth toric surface 
(see \cite[\S 2.5]{fulton})
and $\hat{S}$ its blowup at a torus fixed point.
The exceptional curve of $\psi:\hat{S}\to S$
is denoted by $E$.
Let $X$ be the total space of 
the canonical bundle $K_S$ of $S$
and $\hat{X}=K_{\hat{S}}$.
These are TCY threefolds 
and $E$ is a $(-1,-1)$-curve in $K_{\hat{S}}$. 

Since all torus invariant curves in $X$
are contained in $S\subset X$,
there is a canonical map 
$N_1^T(X)\to H_2(S,\mathbb{Z})$.
In fact, 
the following maps $p,\hat{p}$ are isomorphisms:
\begin{equation}\label{eq:p}
p:H_2(S,\mathbb{R})
\stackrel{\sim}{\longrightarrow}
N_1^T(X)\otimes \mathbb{R},
\qquad
\hat{p}:
 H_2(\hat{S},\mathbb{R})
\stackrel{\sim}{\longrightarrow}
N_1^T(\hat{X})\otimes \mathbb{R}.
\end{equation}

\begin{proposition}\label{prop:blowup}
\begin{enumerate}

\item[(i)] 
For $\beta\in H_2(\hat{S},\mathbb{Z})$ 
such that 
$\beta$ is not an multiple of $[E]$
and satisfying
$\beta.E<0$, 
\begin{equation}\notag
N_{g,\hat{p}(\beta)}(\hat{X})=0~.
\end{equation}

\item[(ii)] 
For $\beta\in H_2(\hat{S},\mathbb{Z})$ such that
$\beta.E=0$,
\begin{equation}\notag
N_{g,\hat{p}(\beta)}(\hat{X})=N_{g,p(\psi_*(\beta))}(X)~.
\end{equation}

\item[(iii)] For a multiple of $[E]$, 
\begin{equation}\notag
N_{g,d[E]}(\hat{X})=
N_{g,d[\mathbb{P}^1]}(\mathcal{O}_{\mathbb{P}^1}(-1)
 \oplus\mathcal{O}_{\mathbb{P}^1}(-1))
.
\end{equation}

\end{enumerate}
\end{proposition}

\begin{proof}

Let $\bar{X}$ be the TCY threefold obtained from $\hat{X}$
by flopping the curve $E$.
Let $\Sigma,\bar{\Sigma},\hat{\Sigma}$ be 
fans of $X,\bar{X},\hat{X}$,
and let $\hat{\tau}_0$  be the 2-cone in $\hat{\Sigma}$
representing $E$.
Then near $\hat{\tau}_0$,
$\hat{\Sigma}$ looks like the right diagram in Figure 
\ref{fig:blowup-fan}
and $\Sigma,\hat{\Sigma}$ are like 
the left and the middle diagrams.
$\Sigma,\bar{\Sigma},\hat{\Sigma}$ 
are identical outside these parts.

A natural inclusion $\Sigma\hookrightarrow\bar{\Sigma}$
induces the isomorphism
\begin{equation}\notag
\alpha:N_1^T(X)\stackrel{\sim}{\longrightarrow}
\{Z\in N_1^T(\bar{X})\,|\,Z.D_{\rho_4}=0\}~,
\end{equation}
where $\rho_4$ is the 1-cone in $\bar{\Sigma}$ shown in 
Figure \ref{fig:blowup-fan}.
Composed with the isomorphism
$\phi_*:N_1^T(\bar{X})\otimes \mathbb{R}
\to N_1^T(\hat{X})\otimes \mathbb{R}$
induced from the flop $\phi:\bar{X} \dashrightarrow  \hat{X}$,
\begin{equation}\label{eq:phi-alpha}
\phi_*\circ\alpha:
N_1^T(X)\otimes \mathbb{R}
\stackrel{\sim}{\longrightarrow}
\{Z\in N_1^T(\hat{X})\otimes \mathbb{R}\,|\,Z.D_{\rho_4}=0\}~,
\end{equation}
where $\rho_4$ is the 1-cone in $\hat{\Sigma}$ shown in 
Figure \ref{fig:blowup-fan}.
By calculating intersection numbers, we see that
the RHS is spanned by (recall $E=C_{\hat{\tau}_0}$)
\begin{equation}\notag
[E]+[C_{\tau_1}],\qquad
[E]+[C_{\tau_2}],\qquad
[C_{\tau}]\quad
(\tau\in\hat{\Sigma}_2'\setminus\{\tau_1,\tau_2\}).
\end{equation}
Under the isomorphisms (\ref{eq:p}),
the inverse of the isomorphism (\ref{eq:phi-alpha})
becomes
\begin{equation}
\psi_*:
\{\beta\in H_2(\hat{S},\mathbb{R})\,|\,\beta.E=0\}
\stackrel{\sim}{\longrightarrow}
H_2(S,\mathbb{R}).
\end{equation}

By applying Theorem \ref{prop:flop2}
to $\hat{X}$ and $\bar{X}$,
we obtain (i).
The statement (iii) follows from the second statement of 
Corollary \ref{cor:GW-flop}. 
The first statement of Corollary \ref{cor:GW-flop},
together with the following, implies (ii):
\begin{equation}\notag
N_{g,\beta}(X)=N_{g,\alpha(\beta)}(\bar{X})
\qquad (\beta\in N_1^T(X)),
\end{equation}
by the
construction of partition function (\ref{eq:partition-fct}).
\end{proof}

\appendix
\section{Combinatorial formulae}\label{ap}
We collect some combinatorial formulae which are used in this paper.
Our basic references are \cite{Mac, EK1, ORV, Zhou}.

$\kappa(\m)$ is always even and 
\beq
\kappa(\m^t) = -\kappa(\m)~, 
\label{kappa}\eeq
where $\mu^t$ denotes the conjugate partition 
(the partition corresponding to the transposed Young diagram of $\m$).

$C_k(\m, \n)$ are nonnegative integers 
which are nonzero for finitely many values of $k$
\cite[Theorem 5.1]{Zhou},
and have the following properties (\cite[\S 3.1]{EK1}, \cite[\S 5.3]{Zhou}):
\begin{equation}
\sum_k C_k(\m, \n) = |\m| +|\n|~, ~\hskip4mm 
\sum_k kC_k(\m,\n)=\frac{1}{2}(\kappa(\m)+ \kappa(\n))~,\\
\label{eq:sumc}
\end{equation}
\begin{equation}
C_k(\m, \n) =C_{-k}(\m^t, \n^t)~.
\label{eq:c}
\end{equation}
The following lemma is proved in
\cite[Lemma in \S C]{EK1}, \cite[Proposition 6.1]{Zhou}:
\begin{lemma}\label{lem:h-C}
For $\m,\n\in\mathcal{P}$,
the following identity holds:
\begin{align*}
\prod_{i,j \geq 1} \left( 1 - Q q^{h_{\m,\n}(i,j)}\right)
&= Z_{(-1,-1)}(q,Q)
\prod_k \left( 1 - Q q^k\right)^{C_k(\m, \n)}.  \label{lemma}
\end{align*}
\end{lemma}

Here are  some  properties of skew Schur function.
Let $x = (x_1, x_2, \dots)$ and $y = (y_1, y_2, \dots)$ be 
sets of variables and $(x, y) = (x_1, x_2, \dots, y_1, y_2, \dots)$.
The following formulae are useful in 
performing the summations over partitions (\cite[p.93,(5.10)]{Mac}):
\beqa
\sum_{\lam\in \mathcal{P}} s_{\lam/\lam_1}(x) s_{\lam/\lam_2}(y) &=&
 \prod_{i,j \geq 1} (1- x_i y_j)^{-1}
\sum_{\m\in\mathcal{P}} s_{\lam_2/\m}(x) s_{\lam_1/\m}(y)~, \label{Schur1} \\
\sum_{\lam\in\mathcal{P}} s_{\lam/\lam_1}(x) s_{\lam^t/\lam_2}(y) &=& 
\prod_{i,j \geq 1} (1+ x_i y_j) 
\sum_{\m\in\mathcal{P}} s_{\lam_2^t/\m}(x) s_{\lam_1^t/\m^t}(y)~, 
\label{Schur2}
\eeqa
\begin{eqnarray} 
&& \sum_{\xi\in\mathcal{P}} s_{\mu/\xi}(x) s_{\xi/\nu}(y) = s_{\mu/\nu}(x, y)~.
\label{Schur3}\end{eqnarray}
Other properties are (\cite[Proposition 4.1]{Zhou}):
\beq
s_{\mu/\nu}(Qx) = Q^{|\mu|-|\nu|}s_{\mu/\nu}(x)~,
\label{homogeneity}\eeq
where $Qx = (Qx_1, Qx_2, \dots)$.
\beq
s_{\lambda/\mu}(q^{\nu+\rho})
 =  (-1)^{|\lambda|-|\mu|}s_{\lambda^t/\mu^t}(q^{-\nu^t-\rho})~.
\label{symmetry}
\eeq

%

\end{document}